\input amstex
\documentstyle{amsppt}


\nopagenumbers

\magnification=\magstep1

\def\alp{\alpha}                
\def\bet{\beta}         
\def\gam{\gamma}                
                \def\Del{\Delta}
\def\eps{\varepsilon}

                \def\Ome{\Omega}


\def\calD{{\Cal D}}

\def\calF{{\Cal F}}



\def\RR{\Bbb R}

\def\CC{\Bbb C}

\def\obet{\overline\bet}
\def\tilrho{{\widetilde{\rho}}}
\def\tilP{{\widetilde{P}}}

\def\sdp{\times \hskip -0.3em {\raise 0.3ex
\hbox{$\scriptscriptstyle |$}}} 


\def\End{\operatorname{End}}
\def\dist{\operatorname{dist}}
\def\olim{\operatornamewithlimits{\overline{lim}\,}}


\def\we{\bigwedge^\b(T^*M)}
\def\wec{\Ome^\b_c(M)}
\def\lwe{L^2\Ome^\b(M)}

\def\b{\bullet}

\def\>{\rangle}
\def\<{\langle}

\NoBlackBoxes
\document


\topmatter

\title On self-adjointness of a Schr\"odinger operator on differential
	forms 
\endtitle

\rightheadtext{On self-adjointness of a Schr\"odinger operator}
\author  Maxim Braverman  \endauthor 

\address
School of Mathematical Sciences,
Tel-Aviv University,
Ramat-Aviv 69978, Israel
\endaddress
\email 
	maxim\@math.tau.ac.il
\endemail

\thanks{The research was  supported by US - Israel Binational Science
Foundation grant No. 9400299}
\endthanks
  \endsubjclass

\abstract
Let $M$ be a complete Riemannian manifold and let $\Ome^\b(M)$ denote
the space of differential forms on $M$. Let $d:\Ome^\b(M)\to
\Ome^{\b+1}(M)$ be the exterior differential operator and let
$\Del=dd^*+d^*d$ be the Laplacian. We establish a sufficient condition
for the Schr\"odinger operator $H=\Del+V(x)$ (where the potential
$V(x):\Ome^\b(M)\to \Ome^\b(M)$ is a zero order differential operator)
to be self-adjoint. Our result generalizes a theorem by I.~Oleinik
about self-adjointness of a Schr\"odinger operator which acts on the
space of scalar valued functions. 
\endabstract
\endtopmatter

\subheading{1. Introduction} Suppose $M$ is a complete Riemannian
non-compact manifold. We will assume that $M$ is oriented and
connected.  Let $T^*M$ denote the cotangent bundle to $M$ and let $\we
= \bigoplus_i \bigwedge^i(T^*M)$ denote the exterior algebra of $T^*M$. We
denote by $\lwe$ the space of square integrable complex valued
differential forms on $M$, i.e. the space of sections of $\we\otimes
\CC$ which are square integrable with respect to the scalar product %
$$
	\<\alp,\bet\> \ = \ \int_M \, \alp\wedge *\obet,
			\qquad \alp,\bet\in \lwe.
\tag1 $$
Here $*$ denotes the Hodge operator associated to the Riemannian
metric on $M$. Note that $L^2\Ome^0(M)$ is just the space of square
integrable complex valued functions on $M$.

Let $d:\lwe\to L^2\Ome^{\b+1}(M)$ denote the exterior differential and
let $d^*$ be the operator formally adjoint to $d$ with respect to the
scalar product (1).

Let $\Del=dd^*+d^*d$ be the Laplacian and consider the Schr\"odinger
operator 
$$
	H \ = \ \Del \ + \ V(x): \ \lwe \ \to  \ \lwe
\tag2  $$
where the {\it potential} $V(x)$ is a measurable section of the bundle
$\End\big(\we\big)$ of endomorphisms of $\we$ which belongs to the
class $L^\infty_{loc}$ (i.e. such that for any compact set $K\subset M$
there exists a constant $C_K>0$ such that $|V(x)|\le C_K$ for almost all 
$x\in K$). 

We denote by $H_0$ the restriction of $H$ on the space $\wec$ of
smooth differential forms with compact support.  The purpose of this
paper is to introduce a sufficient condition on the potential $V(x)$
for operator $H_0$ to be self-adjoint.

\subheading{2. Statement of results} For $x,y\in M$ let $\dist(x,y)$
denote the Riemannian distance between $x$ and $y$.  Fix a point $p\in
M$ and set $r(x)=\dist(x,p)$.
 
Fix $x\in M$. The Riemannian metric on $M$ defines a scalar product
$\<\cdot,\cdot\>_x$ on the fiber $\bigwedge^\b(T^*_xM)\otimes \CC$ of
the vector bundle $\we\otimes \CC$. As usual, we write $V(x)\ge C$ if
$$
	\<V(x)\, \xi,\xi\>_x \ \ge \  C\, \<\xi,\xi\>_x
\tag3  $$
for any $\xi\in \bigwedge^\b(T^*_xM)\otimes \CC$. Note that it follows
from  (3) that $V(x)$ is a self-adjoint endomorphism of
$\bigwedge^\b(T^*_xM)$.

\proclaim{Theorem A} Assume that for almost all $x\in M$ the potential
   $V(x)$ of the operator (2) satisfies the estimate
   $$
	 V(x) \ \ge \  -Q(x),
   \tag4  $$                
   where $1\le Q(x)\le \infty$ and $Q^{-1/2}(x)$ is a Lipschitz function
   on $M$ such that 
   $$
	|\, Q^{-1/2}(x) \ - \ Q^{-1/2}(y)\, | \ \le  \ K\dist(x,y) 
		\quad \text{for any} \quad x,y\in M.
   \tag5  $$
   If for any piecewise smooth curve $\gam:[0,\infty)\to M$ such that
   $\lim_{t\to\infty} r(\gam(t))=\infty$ the integral
   $$
	 \int_\gam   \, 
		Q^{-1/2}(x) \, d\gam \ = \ \infty
   \tag6  $$ 
   then the operator $H_0$ is essentially self-adjoint.  
\endproclaim
For the case of a Schr\"odinger operator acting on scalar valued
functions this theorem was established by I.~Oleinik \cite{O2}.
Note that $Q(x)$ may be equal to infinity on a set of positive
measure.  

As a simple consequence of Theorem~A we obtain the following

\proclaim{Theorem B} Suppose that for almost all $x\in M$ the
   potential $V(x)$ satisfies the estimate $V(x)\ge -q\big(r(x)\big)$,
   where $1\le q\le \infty$ and $q^{-1/2}(t)$ is a Lipschitz function
   on $\RR$ such that $\int_0^\infty \, q^{-1/2}(t) \, dt = \infty$.
   Then the operator $H_0$ is essentially self-adjoint.

  In particular, if $M=\RR^n$ and $V(x)\ge -C|x|^2$ then the operator 
  $H_0$ is essentially self-adjoint.
\endproclaim

\subheading{\it Remark} Theorem~A remains true if we replace
$\lwe$ by the space of square integrable forms on $M$ with values in a
flat Hermitian vector bundle $\calF$ over $M$, provided that the
Hermitian structure on $\calF$ is flat. In this case the differential
$d$ should be replaced by the covariant differential associated to the
flat structure on $\calF$. The proof is a verbatim repetition of the
proof for the scalar case, cf. bellow. However the notation in the
vector valued case is more complicated.

\subheading{3. Historical remarks} An analogue of Theorem~B for the
case \/ $M=\RR^1$ \/ was established by Sears \cite{Se}.  B.~Levitan
\cite{Le} proved the Sears theorem for the Schr\"odinger operator
acting on scalar valued functions on $M=\RR^n$.  F.~Rofe-Beketov
\cite{RB} extended these results to the case where the potential
$V(x)$ can not be estimate by a function depending only on
$\dist(x,p)$. Many results and references about the essential
self-adjointness of Schr\"odinger operators on $\RR^n$ may be found in
\cite{RS}.

I.~Oleinik \cite{O1,O2} established Theorem~A for the Schr\"odinger
operator acting on scalar valued functions on a complete Riemannian
manifold.

Essential self-adjointness of a pure Laplacian (without lower order
terms) on differential forms on a complete Riemannian manifold was
first stated and proved by M.~P.~Gaffney \cite{Ga1,Ga2}. A number of
related results may be found in \cite{Sh}.

In \cite{BFS}, Theorem~B is established for the case where $M$ is
a manifold with cylindrical ends and the potential $V(x)\ge 0$. The
result is used there to study Witten deformation of the Laplacian on a
non-compact manifold.

\subheading{4. Acknowledgment} I would like to thank  M.~Shubin for
posing the problem and for essential help. 

I am very grateful to I.~Oleinik for pointing out a gap in a
preliminary version of the paper and for drawing my attention to the 
paper \cite{O2}.

I am also thankful to M.~Farber and V.~Matsaev for valuable
discussions.

\subheading{5. The domain of $\calD(H_0^*)$} Let $H_0^*$ denote the
operator adjoin to $H_0$. The domain $\calD(H_0^*)$ of $H_0^*$
consists of forms $\alp\in \lwe$ such that $H\alp$ understood in the
sense of distributions also belongs to $\lwe$.

The operator $H_0$ is symmetric. Hence, to show that its closure is
self-adjoint it is enough to show that the adjoint operator $H_0^*$ is
symmetric. In other words we have to prove that
$$ 
	\int_M \, 
	 (H\alp\wedge *\obet   -  \alp\wedge *\overline{H\bet})  
    \ = \ 0
\tag7  $$
for any $\alp, \bet\in \calD(H_0^*)$.

To prove (7) we need some information about the behavior of
differential forms from $\calD(H_0^*)$.  The main result of this
section is the following lemma, which provides us with this
information.
\proclaim{Lemma 1} If $\alp\in \calD(H_0^*)$ then the forms
   $Q^{-1/2}d\alp, Q^{-1/2}d^*\alp$ are square integrable. 
\endproclaim
\subheading{\it Remark} 1. \ By the standard theory of elliptic
operators any $\alp\in \calD(H_0^*)$ belongs to the Sobolev space
$H_{loc}^2$. Hence, $d\alp, d^*\alp$ are locally square
integrable. Thus the lemma provides us with an information about the
behavior of the forms from $\calD(H_0^*)$ at infinity.

2. \ For the Schr\"odinger operator on scalar valued functions on
$\RR^n$ an analogous lemma was established in \cite{RB}. The proof was
adopted in \cite{O1,O2} to the case of a Riemannian manifold. In our
proof we follow rather closely the lines of \cite{O2}. However, the
fact that we deal with differential forms rather than with scalar
valued functions demands a more careful analysis.

\demo{Proof}
Recall that we fixed a point $p\in M$ and that for any $x\in M$ we
denoted by $r(x)$ the Riemannian distance between $x$ and $p$. 
 
It is shown in \cite{O2, Proof of Lemma~1} that for any $R>0,\eps>0$
there exist smooth functions $r_{R,\eps}(x),F_{R,\eps}(x)$ on $M$
which approximate the Lipschitz functions $r(x),Q^{-1/2}(x)$ in the
sense that 
$$
   \aligned
	|r_{R,\eps}(x)-r(x)| \ < \ \eps, \qquad
		Q^{-1/2}(x)-\eps \ &< \ 
		    F_{R,\eps}(x) \ < \ (1+\eps) \, Q^{-1/2}(x),\\ 
	 \olim_{\eps\to 0} |d r_{R,\eps}(x)|\le 1, \qquad \ \ 
		&\olim_{\eps\to 0} |dF_{R,\eps}(x)| \ \le \ K,
   \endaligned
\tag8  $$
for any  \/ $x\in r_{R,\eps}^{-1}([0,R+1])$. Here $K$ is the same
constant as in (5).

Let $\Psi:[0,+\infty)\to [0,1]$ be a smooth function which is equal
to one when $t\le 1/2$ and  which is equal to zero when $t\ge 1$. Set
$$
	\psi_{R,\eps}(x) \ = \ 
	   \cases
		\Psi\big(\frac{r_{R,\eps}(x)}R\big)\, F_{R,\eps}(x)
		  \quad &\text{if} \quad r_{R,\eps}(x)\le R; \\
		0 \quad &\text{outside of the set} \  r_{R,\eps}(x)\le R.
	   \endcases
\tag9  $$

For any $R>0$ the functions $\psi_{R,\eps}, \, \eps<1$ vanish outside
of the compact set $r^{-1}([0,R+1])$. Hence, it follows from (8) and
(4) that there exist a constant $K_1>0$ not depending on $R$ and a
number $\eps_R>0$ (which does depend on $R$) such that
$$
	 |d \psi_{R,\eps}(x)| \ \le \  K_1, \qquad 
		\psi_{R,\eps}^2(x) \ \le \ 2, \qquad
	   \left|\, \int_M\, \psi_{R,\eps}^2\, \alp\wedge *V\alp\, \right| 
			\ \le \   2\, \|\alp\|^2,
\tag10  $$
for any $x\in M, \, R>1, \, 0<\eps<\eps_R, \, \alp\in \lwe$.  Here
$\|\alp\|=\<\alp,\alp\>^{\frac12}$ denotes the $L^2$-norm of the form
$\alp$.

Functions \/ $\psi_{R,\eps}$ \/ have compact support. Hence, in view
of the remark~1 after the statement of the lemma, the forms \/
$\psi_{R,\eps}d\alp$ \/ and \/ $\psi_{R,\eps}d^*\alp$ \/ are square
integrable. Assume that \/ $\alp\in \calD(H_0^*)$ \/ is a real valued
form and set
$$
	  J_{R,\eps}^2 \ = \ 
	 \|\psi_{R,\eps}d\alp\|^2  +  \|\psi_{R,\eps}d^*\alp\|^2 \ = \
		\int_M \psi_{R,\eps}(x)^2\, 
	   \big(d\alp\wedge *d\alp + d^*\alp\wedge *d^*\alp\big).
\tag11  $$

It follows from (8), (9) that to prove the lemma it is
enough to show that
$$
	\olim_{R\to \infty} \, \olim_{\eps\to 0} J_{R,\eps} 
			\ < \ \infty.
\tag12  $$

Let us first rewrite the integrand in (11) in a more convenient
form. In the calculations bellow we  use the equality (cf.
\cite{Wa, \S6.1}) \/ $d^*\alp= (-1)^{|\alp|}*^{-1}d*\alp$ \/ where
$|\alp|$ denotes the degree of the differential form $\alp$. 
$$
  \alignat 2
    \  \psi_{R,\eps}^2\, d\alp\wedge &&*d\alp  \ =& \
         d\big(\psi_{R,\eps}^2\alp\wedge*d\alp\big) - 
		2\psi_{R,\eps}d\psi_{R,\eps}\wedge\alp\wedge *d\alp +
			\psi_{R,\eps}^2\alp\wedge*d^*d\alp,  \\
     \  \psi_{R,\eps}^2\, d^*\alp\wedge &&*d^*\alp \  =& \ 
	 (-1)^{|\alp|} \psi_{R,\eps}^2\, d^*\alp\wedge d*\alp  \tag13 \\ 
	 &&=   -d&\big(\psi_{R,\eps}^2d^*\alp\wedge*\alp\big)  + 
		2 \psi_{R,\eps}d\psi_{R,\eps}\wedge d^*\alp\wedge *\alp +
			\psi_{R,\eps}^2dd^*\alp\wedge *\alp.  
  \endalignat
$$
It follows now from (13), (10) and from the Stokes theorem that, if
$R>1, \, \eps<\eps_R$, then
$$
  \alignat 2
	\|\psi_{R,\eps}d\alp\|^2 
	\  &=&& \ 
		\int_M \,  \psi_{R,\eps}^2\, d\alp\wedge *d\alp 
	\ = \
		-2\<d\psi_{R,\eps}\wedge \alp,\psi_{R,\eps}d\alp\> + 
			\<\alp,\psi_{R,\eps}^2d^*d\alp\>   
	  \\
	&\le&& \ 
		2K_1\|\alp\|\, \|\psi_{R,\eps}d\alp\| + 
			\<\alp,\psi_{R,\eps}^2d^*d\alp\>,     
	  \\
	 \|\psi_{R,\eps}d^*\alp\|^2 
	\  &=&& \ 
		\int_M \, \psi_{R,\eps}^2\, d^*\alp\wedge *d^*\alp 
	\ = \
		   2\<d\psi_{R,\eps}\wedge\psi_{R,\eps}d^*\alp,\alp\> + 
			\<\psi_{R,\eps}^2dd^*\alp,\alp\>  
		  \\
        &\le&& \
		2K_1 \|\psi_{R,\eps}d^*\alp\|\, \|\alp\| +
				 \<\alp,\psi_{R,\eps}^2dd^*\alp\>.
  \endalignat
$$
Summing these two equations we obtain
$$
  \multline
	J_{R,\eps}^2 \ \le \ 
	  2K_1\|\alp\| \big(
		\|\psi_{R,\eps}d\alp\| + \|\psi_{R,\eps}d^*\alp\|
								\big) +
		\<\alp,\psi_{R,\eps}^2\Del\alp\> \\ 
      \le \
	  4K_1\|\alp\|\,  J_{R,\eps} + 
	     \int_M\, 
	     \psi_{R,\eps}^2\big(\alp\wedge *H\alp - \alp\wedge *V\alp\big) \\
     \le \
	4K_1\|\alp\|\,  J_{R,\eps} + 
		2 \|\alp\|\, \|H\alp\| +2 \|\alp\|^2.
  \endmultline
\tag14 $$
Here the last inequality follows from (10).

It follows from (14) that the set $\{J_{R,\eps}: \, R>1,
\eps<\eps_R\}$ is bounded from above. Hence (12) holds. The proof of
the lemma is completed.  
$\square$\enddemo

\subheading{6. Proof of Theorem~A} We apply a modification of the
method used in \cite{RB} suggested by I.~Oleinik \cite{O2}.

The quantity
$$
	\tilrho(x,y) = \inf_{\gam}\, \int_\gam Q^{-1/2}(x)\, d\gam,
\tag15  $$
where the infimum is taken over all piecewise smooth curves connecting
the points $x,y\in M$, is called {\it generalized distance} between $x$
and $y$. It is a symmetric function in $x,y$ which satisfies the
triangular inequality. The first metric axiom is not valid in
general. Note, however, that (6) implies, that the sets
$P^{-1}([0,R])$ are compact for any $R>0$. 

Recall that in Section~2 we have fixed a point $p\in M$.
Set \/ $P(x)=\tilrho(x,p)$. Then (cf. \cite{O2, Lemma~2})
$$
	|P(x)-P(y)| \ \le \ 
		Q^{-1/2}(x)\dist(x,y) \ + \frac{K}2 (\dist(x,y))^2
\tag16  $$
for any $x,y\in M$. It follows (cf. \cite{O2}) that for any $R>0,
\eps>0$ there exists a smooth function $\tilP_{R,\eps}(x)$ which
approximates $P(x)$ in the sense that
$$
	|\tilP_{R,\eps}(x)-P(x)|\le \eps, \qquad 
	  \olim_{\eps\to 0} |d\tilP_{R,\eps}(x)|\le Q^{-1/2}(x),
\tag17  $$
for any $x\in P^{-1}([0,R+1])$. 

Assume that $\eps<1$ so that $\tilP_{R,\eps}^{-1}([0,R])\subset
P^{-1}([0,R+1])$. Let us define a piecewise smooth function
$P_{R,\eps}(x)$ on $M$ by the formula 
$$
	P_{R,\eps}(x)=\cases
		\tilP_{R,\eps}(x) \quad &\text{if} \quad 
					\tilP_{R,\eps}(x)\le R; \\
		R	   \quad &\text{outside the set} \quad
					 \tilP_{R,\eps}(x)\le R.
		\endcases
\tag18  $$
By (17), the inequality
$$
	\olim_{\eps\to 0} |dP_{R,\eps}(x)| \ \le \ Q^{-1/2}(x)
\tag19  $$
holds almost everywhere on $M$. 

Recall from Section~5, that the statement of Theorem~A is equivalent
to equality (7).  Fix $\alp,\bet\in \calD(H_0^*)$ and consider the
following approximation of the integral (7)
$$\multline
	I_{R,\eps} \ = \ \int_M \, \left(1-\frac{P_{R,\eps}}R\right)\, 
		\big(H\alp\wedge*\obet -
			\alp\wedge*\overline{H\bet}\big) \\
	\ = \
	\int_M \, \left(1-\frac{P_{R,\eps}}R\right)\, 
		\big(\Del\alp\wedge*\obet -
			\alp\wedge*\overline{\Del\bet}\big).
  \endmultline
\tag20  $$
By the Fatou theorem (\cite{RS, Theorem I.17}), it is enough to show that 
$$
	\olim_{R\to \infty} \, \olim_{\eps\to 0} I_{R,\eps} \ = \ 0.
\tag21  $$

We will need the following ``integration by parts'' lemma
\footnote{I learned this lemma from M.~Shubin.}
\proclaim{Lemma 2} Let $\phi:M\to \RR$ be a smooth function with compact
   support.  Then 
   $$
    \multline
	\ \int_M\, \phi \, \Del\alp\wedge *\obet  \\
      \ = \
	\int_M\, 
	  \phi \, \big(d\alp\wedge*d\obet + d^*\alp\wedge*d^*\obet\big)
	\ + \
	      \int_{M} \, d\phi\wedge \big(
		  \obet\wedge *d\alp -  d^*\alp\wedge*\obet \big)
     \endmultline \tag22
   $$	
   for any $\alp,\bet\in \calD(H_0^*)$.
\endproclaim
Note that, by remark 1 after the statement of Lemma~1, all the
integrals in (22) have sense.
\demo{Proof} 
Recall that \/ $d^*u= (-1)^{|u|}*^{-1}d*u$ \/ where $|u|$ denotes
the degree of the differential form $u$. Hence, if $|u|=|v|-1$, then    
$$
	\phi du\wedge*v \ = \ \phi u\wedge*d^*v \ - \
		d\phi\wedge u\wedge *w \ + \ d\,(\phi u\wedge*v)
\tag23  $$ 
Substituting into (23) first \/ $u=d^*\alp,\, v=\obet$ \/ and
then \/ $u=\obet,\, v=d\alp$ \/ we obtain
$$
  \align
	\phi dd^*\alp\wedge*\obet \ &= \ 
	  -d\phi\wedge d^*\alp\wedge *\obet \ + \
	    \phi d^*\alp\wedge*d^*\obet \ + \ 
		d\, (\phi d^*\alp\wedge*\obet), \\
	\phi d^*d\alp\wedge *\obet \ &= \ \phi \obet\wedge*d^*d\alp \ = \ 
	    d\phi\wedge \obet\wedge d\alp \ + \
		\phi d\alp\wedge *d\obet \ - \ 
				d\, (\phi \obet\wedge*d\alp).
  \endalign
$$
In the last equality we used that $u\wedge*v=v\wedge *u$ for any
differential forms $u,v$ of the same degree. Summing the above
equations, integrating over $M$  and using the Stokes theorem we get 
(22).
$\square$\enddemo

Using  definition (20) of $I_{R,\eps}$ and Lemma~2 we
obtain
$$
	I_{R,\eps} \ = \   \frac1R\, \int_M \,	dP_{R,\eps}\wedge \big(
		  \obet\wedge *d\alp -  d^*\alp\wedge*\obet -
		\alp\wedge *d\obet + d^*\obet\wedge*\alp \big).
\tag24  $$

Let \/ $d\mu(x)$ \/ denote the Riemannian density on \/ $M$. For any
\/ $\xi\in \bigwedge^k(T^*M)\otimes \CC$ \/ we denote by \/ $|\xi|$ \/
its norm with respect to the scalar product on \/
$\bigwedge^\b(T^*M)\otimes\CC$ \/ induced by the Riemannian structure
on $M$. Then
$$
	|\<\alp,\bet\>| \ \le \ \int_M\, |\alp\wedge*\obet|\, d\mu(x) 
	\ \le \ \int_M\, |\alp|\, |\bet|\, d\mu(x) \ \le \ 
	\|\alp\|\, \|\bet\|
\tag25  $$
for any $\alp,\bet\in \lwe$. 

Let us estimate the behavior of the right hand side of (24)
as $\eps\to 0$. For the first term we obtain
$$
  \multline
	\olim_{\eps\to 0}
	\left|\, 
	   \frac1R \int_M\, dP_{R,\eps}\wedge \obet\wedge *d\alp 
	\right|    
	\ \le \
	\frac1R\, \olim_{\eps\to 0} \int_M \, 
		|dP_{R,\eps}|\, |d^*\alp|\, |\obet|\, d\mu(x) 
	\\ \le \
	\frac1R\,  \int_M \, 
		|Q^{-1/2} d^*\alp|\, |\obet|\, d\mu(x) 
	\ \le \
		\frac{\|Q^{-1/2}d^*\alp\|\, \|\bet\|}R.
  \endmultline
\tag26 $$	 
In the second inequality in (26) we used the estimate
(19). The last inequality in (26) follows from
Lemma~1.

Analogously, one can estimate the other terms in the right hand side of 
(24). That proves (21) and Theorem~A. \
$\square$


\enddocument